\documentclass[11pt]{amsart}
\usepackage{amsfonts}
\usepackage{amsmath}
\usepackage{amsthm}
\usepackage{amssymb}
\usepackage{latexsym}
\advance\textwidth by 1.2in \advance\oddsidemargin by -.6in
\advance\evensidemargin by -.6in

\newtheorem{cor}{Corollary}[section]
\newtheorem{lem}{Lemma}[section]

\newtheorem{prop}{Proposition}[section]

\theoremstyle{definition}
\newtheorem{defn}{Definition}[section]
\theoremstyle{definition}
\newtheorem{thm}{Theorem}

\newenvironment{pf}{\proof}{\endproof}
\newcounter{cnt}

\theoremstyle{remark}

\numberwithin{equation}{section} 
\begin{document}

\newcommand{\thmref}[1]{Theorem~\ref{#1}}
\newcommand{\secref}[1]{Sect.~\ref{#1}}
\newcommand{\lemref}[1]{Lemma~\ref{#1}}
\newcommand{\propref}[1]{Proposition~\ref{#1}}
\newcommand{\corref}[1]{Corollary~\ref{#1}}
\newcommand{\remref}[1]{Remark~\ref{#1}}
\newcommand{\er}[1]{(\ref{#1})}
\newcommand{\id}{\operatorname{id}}
\newcommand{\tensor}{\otimes}
\newcommand{\nc}{\newcommand}
\newcommand{\rnc}{\renewcommand}
\newcommand{\qbinom}[2]{\genfrac[]{0pt}0{#1}{#2}}
\nc{\cal}{\mathcal} \nc{\goth}{\mathfrak} \rnc{\bold}{\mathbf}
\renewcommand{\frak}{\mathfrak}
\renewcommand{\Bbb}{\mathbb}
\nc\bpi{{\mbox{\boldmath $\pi$}}} \nc\bolambda{{\mbox{\boldmath
$\lambda$}}} \nc\bomu{{\mbox{\boldmath $\mu$}}}

\newcommand{\lie}[1]{\mathfrak{#1}}

\nc{\Cal}{\cal} \nc{\Xp}[1]{X^+(#1)} \nc{\Xm}[1]{X^-(#1)}
\nc{\on}{\operatorname} \nc{\ch}{\mbox{ch}} \nc{\Z}{{\bold Z}}
\nc{\J}{{\cal J}} \nc{\C}{{\bold C}} \nc{\Q}{{\bold Q}}
\renewcommand{\P}{{\cal P}}
\nc{\N}{{\Bbb N}} \nc\boa{\bold a} \nc\bob{\bold b} \nc\boc{\bold
c} \nc\bod{\bold d} \nc\boe{\bold e} \nc\bof{\bold f}
\nc\bog{\bold g} \nc\boh{\bold h} \nc\boi{\bold i} \nc\boj{\bold
j} \nc\bok{\bold k} \nc\bol{\bold l} \nc\bom{\bold m}
\nc\bon{\bold n} \nc\boo{\bold o} \nc\bop{\bold p} \nc\boq{\bold
q} \nc\bor{\bold r} \nc\bos{\bold s} \nc\bou{\bold u}
\nc\bov{\bold v} \nc\bow{\bold w} \nc\boz{\bold z}

\nc\ba{\bold A} \nc\bb{\bold B} \nc\bc{\bold C} \nc\bd{\bold D}
\nc\be{\bold E} \nc\bg{\bold G} \nc\bh{\bold H} \nc\bi{\bold I}
\nc\bj{\bold J} \nc\bk{\bold K} \nc\bl{\bold L} \nc\bm{\bold M}
\nc\bn{\bold N} \nc\bo{\bold O} \nc\bp{\bold P} \nc\bq{\bold Q}
\nc\br{\bold R} \nc\bs{\bold S} \nc\bt{\bold T} \nc\bu{\bold U}
\nc\bv{\bold V} \nc\bw{\bold W} \nc\bz{\bold Z} \nc\bx{\bold X}

\nc\gfin{\mbox{${\frak g}_{fin}$}} \nc\gaff{\mbox{${\frak
g}_{aff}$}}
 \nc\gtor{\mbox{${\frak g}_{tor}$}}
 \nc\gcurr{\mbox{${\frak g}_{tor}^+$}}

 \nc\getor{\mbox{${\frak g}^e_{tor}$}}

 \nc\hetor{\mbox{${\frak h}^e_{tor}$}}

 \nc\geaff{\mbox{${\frak
g}^e_{aff}$}}

\nc\heaff{\mbox{${\frak h}^e_{aff}$}}

\nc\hfin{\mbox{${\frak h}_{fin}$}} \nc\haff{\mbox{${\frak
h}_{aff}$}}

\nc\htor{\mbox{${\frak h}_{tor}$}}

\nc\npfin{\mbox{${\frak n}^+_{fin}$}}

 \nc\npaff{\mbox{${\frak n}^+_{aff}$}}

\nc\nptor{\mbox{${\frak n}^+_{tor}$}}

\nc\nmfin{\mbox{${\frak n}^-_{fin}$}}

\nc\nmaff{\mbox{${\frak n}^-_{aff}$}}

\nc\nmtor{\mbox{${\frak n}^-_{tor}$}}

\nc\bufin{\mbox{${\bu}_{fin}$}} \nc\buaff{\mbox{${\bu}_{aff}$}}

\nc\bug{\mbox{${\bu}(\frak g_{tor}(>)$}}
\nc\bul{\mbox{${\bu}(\frak g_{tor}(<))$}}
\nc\buz{\mbox{${\bu}(\frak g_{tor}(0))$}}

 \nc\butor{\mbox{${\bu}_{tor}$}}
  \nc\buetor{\mbox{${\bu}^e_{tor}$}}
   \nc\bueaff{\mbox{${\bu}^e_{{\text{aff}}}$}}

\title{Representations of double affine lie algebras}\author{Vyjayanthi Chari and
Thang Le}
\address{ Department of Mathematics, University of
California, Riverside, CA 92521.} \maketitle

\section{Introduction} The representation theory of Kac--Moody
algebras, and in particular that of affine Lie algebras, has
been extensively studied over the past twenty years. The
representations that have had the most applications are the
integrable ones, so called because they lift to the corresponding
group.

The affine Lie algebra associated to a finite-dimensional complex
simple Lie algebra $\frak g$ is the universal (one--dimensional)
central extension  of the Lie algebra of polynomial maps from
$\bc^*\to\frak g$. There are essentially two kinds of interesting
integrable representations of this algebra: one is where the where
the center acts as a positive integer, or positive energy
representations; and the other is where the center acts trivially,
or the level zero. Both kinds of representations have interesting
applications: representations of the first kind have connections
with number theory through the Rogers--Ramanujam identities, the
theory of vertex algebras and conformal field theory;
representations of the second kind are connected with the
six--vertex and $XXZ$--model \cite{Mc} and the references therein,
the Kostka polynomials and the fermionic formulas of Kirillov and
Reshetikhin, \cite{KR}. The study of such level zero
representations was begun in \cite{C} and continued in
\cite{CP},\cite{CPnew}. More recently a geometric approach to such
representations was developed in \cite{N} for the corresponding
quantum algebras.

The affine Lie algebras admit an obvious generalization. Namely,
we can consider central extensions of the polynomial maps
$(\bc^*)^\ell\to\frak g$. Not surprisingly, these algebras are a
lot more complicated, for instance the central extension is now
infinite---dimensional. A systematic study of such algebras can be
found  in \cite{ABP} and the representation theory has been
studied in \cite{BB}, \cite{E1}, \cite{E2}, \cite{EMY}. In
general, interesting theories have been found for the quotients of
this algebra by a central ideal of finite--codimension.

One such quotient is the double affine algebra, this algebra is
obtained from the affine algebra in the same way that the affine
algebra is obtained from the finite--dimensional algebra. One can
also define a corresponding quantum object, and representations of
these have been studied in \cite{N}, \cite{Saito}, \cite{VV}.
However, relatively little is still known about the integrable
representations of the double affine algebras.

In this paper we  study representations of the double affine
algebra $\gtor$ when one of the centers acts trivially; this is
also the situation studied in the quantum case mentioned above.
The category of such representations is not semisimple and  our
interest is in indecomposable integrable representations of
$\gtor$ rather than the  irreducible ones.  We are motivated by
considerations coming from the study of quantum affine algebras
\cite{CPweyl} and modular Lie algebras. Thus, we believe that
these indecomposable representations should be the limit as $q\to
1$ of the irreducible representations of the corresponding quantum
algebra. In the case of quantum affine algebras this is in fact a
conjecture which has been checked in many cases; for double affine
algebras, no such conjecture is possible at the moment, since the
notion of $q\to 1$ is not well--defined. However, Theorem
\ref{irred} can be viewed as the classical analog of the result in
\cite{Saito}.

The paper is organized as follows. In section 1, we set up the
notation to be used in the rest of the paper. In section 2, we
recall several result on the representations of the affine
algebra. In section 3, we define a family of indecomposable
representations and identify their irreducible quotients.

Sections 4 and 5 are devoted to the study of the structure of
these representations. Thus we give a sufficient  condition for
the indecomposable modules to be irreducible. We also show  that
in fact the modules are almost always reducible in Section 5. In
the case of affine algebras, this was proved by passing to  the
quantum situation. In this paper, however, we prove it by using
the notion of fusion product representations, an idea introduced
recently by Feigin and Loktev in \cite{FL}.

\section{Preliminaries}

\subsection{}
Let $\gfin$ be a complex finite-dimensional simple Lie algebra of
rank $n$, and let $\hfin$ be a Cartan subalgebra of $\gfin$.
 Fix a set   $\{\alpha_i:1\le i\le n \}$ (resp.
$\{\omega_i:1\le i\le n\}$)
   of simple roots
 (resp.  fundamental weights) of $\gfin$ with respect to $\hfin$.
  Let $R^+_{fin}$ be the corresponding set of positive roots and
  let
  $\theta\in R^+_{fin}$ be the highest root. Given $\alpha\in R^+_{fin}$, let
$\frak g_{fin}^{\pm \alpha}$ be the corresponding root spaces and
fix non--zero elements $x^\pm_\alpha\in \frak g_{fin}^{\pm\alpha}$
and $h_\alpha\in \hfin$ such that the elements $x^\pm_\alpha$,
$h_\alpha$ span a subalgebra isomorphic to $sl_2$. Set
$x^\pm_i=x^\pm_{\alpha_i}$, $h_i=h_{\alpha_i}$. Define subalgebras
$$\frak n_{fin}^\pm =\oplus_{\alpha\in R^+_{fin}}\frak
g_{fin}^{\pm \alpha}.$$  We have,
$$\gfin=\nmfin\oplus\hfin\oplus\npfin.$$ Let  $Q_{fin}^+$ (resp.
$P_{fin}^+$) be  the non-negative root (resp. weight) lattice of
$\gfin$.

\subsection{} The corresponding untwisted affine Lie algebra $\gaff$
is the universal central extension of the Lie algebra of Laurent
polynomial maps from $\bc^*\to\gfin$. Thus $$\gaff=\gfin\otimes
\bc[t_1,t_1^{-1}]\oplus \bc c_1,$$ with bracket given by $$
[xt_1^r, yt_1^m]=[x,y]t_1^{r+m}+ r\delta_{r,-m}<x,y>c_1,\ \
[c_1,xt_1^r]=0,\ \ x,y\in\gfin, \ r,m\in\bz$$ where $<,>$ is the
Killing form of $\gfin$.  Let $d_1$ be the derivation $t_1d/dt_1$
of $\bc[t_1,t_1^{-1}]$; then $d_1$ acts on $\gaff$ in the obvious
way. The extended affine algebra $\frak g_{ aff}^e $ is the
semi--direct product of $\bc d_1$ and $\gaff$, i.e.
$\geaff=\gaff\oplus\bc d_1$ with the Lie bracket between $d_1$ and
$\gaff$ being given by $$[d_1, c_1]=0,\ \ [d_1,xt_1^m]=mxt_1^m,\
\ x\in\gfin, m\in\bz.$$ The algebra $\geaff$ admits a symmetric,
invariant bilinear form, defined as follows: $$ <xt_1^r\ , \
yt_1^m>_{aff}= \delta_{r,-m}<x\ ,y >\  \ \ x,y\in \gfin,\ \
r,m\in\bz$$ and$$<c_1\ ,\ d_1>=1, \ \ <c_1\ ,\ \gaff>\  =\ <d_1
\,\ \gaff>\ =\ 0.$$

 Set $$\heaff=\hfin\oplus\bc
c_1\oplus\bc d_1, \ \ \npaff =\frak g\otimes
t_1\bc[t_1]\oplus\npfin,\ \ \nmaff=\frak g\otimes
t_1^{-1}\bc[t_1^{-1}]\oplus\nmfin.$$ Clearly we can regard $\gfin$
as the algebra of constant maps in $\geaff$. Given any element
$\lambda\in (\hfin)^*$ we regard it as an element of $(\heaff)^*$
by setting $\lambda(c_1) = \lambda(d_1)=0$. Define
$\delta_1\in(\heaff)^*$ by $$\delta_1(d_1)=1, \ \
(\delta_1)|_{\hfin\oplus\bc c_1}=0.$$ Set
$\alpha_0=\delta_1-\theta$. The elements  $\{\alpha_i: 0\le i\le
n\}$ are called the affine simple  roots. The corresponding set of
positive affine roots is $$R^+_{aff}=\{\pm
\alpha+r\delta_1:\alpha\in R_{fin}^+, r\in\bz, r>0\}\cup
R^+_{fin}\cup \{r\delta_1: r\in\bz,r>0\}.$$  The root space
corresponding to $ \pm\alpha+r\delta_1$ is $\frak g_{fin}^{\pm
\alpha}\otimes \bc t_1^r$, and that corresponding to $r\delta_1$
is $\hfin\otimes \bc t_1^r$. Fix non--zero elements
$x_0^\pm\in\frak {g}_{aff}^{\pm\alpha_0}$ and $h_0\in\heaff$ so
that the subalgebra spanned by $x_0^\pm, h_0$ is isomorphic to
$sl_2$. Then,  $$c_1=h_0+h_\theta.$$
 Let $\omega_0\in(\heaff)^*$
be defined by $$\omega_0(h)=0,\ \ \omega_0(c_1)=1, \ \
\omega_0(d_1)=0,\ \ (h\in\hfin).$$ The affine root lattice
$Q^+_{aff}$ and the affine weight lattice $P^+_{aff}$ in
$(\heaff)^*$ are now  defined in the obvious way.

\subsection{} The double affine algebra $\gtor$ and the extended
algebra $\getor$ are  obtained from $\gaff$ in the same way that
$\gaff$ and $\geaff$ were  obtained from $\frak g_{fin}$.  Thus,
$$\getor=\gaff\otimes\bc[t_2,t_2^{-1}]\oplus\bc c_2\oplus\bc
d_1\oplus \bc d_2, \ \ \gtor
=\gaff\otimes\bc[t_2,t_2^{-1}]\oplus\bc c_2\oplus\bc d_1,$$ with
the Lie bracket given as follows. The element $c_2$ is central,
and $$[d_2,xt_2^r]= rxt_2^r,\ \
[xt_2^r,yt_2^m]=[x,y]t_2^{r+m}+\delta_{r,-m}<x,y>_{aff}c_2,\ \ \
(x,y\in\gaff, r,m\in\bz). $$ Clearly  $\gaff$ is a subalgebra of
$\gtor$. Set $$\hetor=\heaff\oplus\bc c_2\oplus \bc d_2.$$ As
before, we regard an element  $\lambda\in (\heaff)^*$ as an
element of $(\hetor)^*$ by setting it to be zero on $\bc c_2$ and
$\bc d_2$. The set of roots of the double affine  algebra are
then, $$ R_{tor}=\{\pm\alpha+n_2\delta_2: \alpha\in R^+_{aff},
n_2\in\bz\}\cup\{n_2\delta_2:n_2\in\bz\}.$$ Unlike the affine
case, there is no natural choice of simple roots. However, we will
work with the following partition of $R_{tor}$ into mutually
disjoint sets,
\begin{eqnarray*}
&R_{tor}(>)&=\{\alpha+n_2\delta_2: \alpha\in R^+_{aff},
n_2\in\bz,\}\\ & R_{tor}(0)&= \{n_2\delta_2:n_2\in\bz\},\\
&R_{tor}(<)&=\{-\alpha+n_2\delta_2: \alpha\in R^+_{aff},
n_2\in\bz,\}.\end{eqnarray*} The subalgebras $\frak g_{tor}(>)$,
$\frak g_{tor}(0)$ and $\frak g_{tor}(<)$ are defined in the
obvious way.

Given a Lie algebra $\frak a$,  let $\bu(\frak a)$ denote its
universal enveloping algebra. We let $\bu_{fin}$ (resp. $\buaff$,
$\butor$) denote the enveloping algebra of $\gfin$ (resp.
$\geaff$, $\gtor$). By the Poincare--Birkhoff--Witt theorem, we
have
\begin{eqnarray*}&\bufin &=\bu(\nmfin)\bu(\hfin)\bu(\npfin),\\
&\buaff&=\bu(\nmaff)\bu(\heaff)\bu(\nmaff),\\
&\butor&=\bu(\frak g_{tor}(<))\bu(\htor)\buz\bu(\frak
g_{tor}(>).\end{eqnarray*}

For $h\in\haff$, let $\Lambda^\pm(h,u)$  be the power series in an
indeterminate  $u$ with coefficients in $\butor$, defined by
\begin{equation}\label{power}  \Lambda^\pm(h,u)=
{\text{exp}}-\left(\sum_{r=1}^\infty \frac{h t_2^{\pm
r}u^r}{r}\right). \ \
\end{equation}
Let $\Lambda^\pm(h,r)$ be the coefficient of $u^{r}$ in
$\Lambda^\pm(h,u)$. It is not hard to see  that the elements
$\{\Lambda^\pm(h, r): h\in \haff, r\in\bz, r>0\} $ generate the
subalgebra $\buz$.

For $\alpha\in R^+_{fin}$ and $r\in \bz$, the
 the elements
$\{x^\pm_\alpha t_1^r  t_2^s:s\in\bz\}$ generate a subalgebra of
$\gtor$ which is isomorphic to  the affine algebra associated to
$sl_2$.
 The following lemma was proved in \cite{G}.
 \begin{lem}\label{deflambda}
  Let $\beta\in R^+_{aff}$ be of the
 form $\alpha+r_1\delta_1$ for some $r_1\in\bz$.
 For all  $s\ge
 1$, we have
 \begin{eqnarray*}&(x^+_{\beta}t_2^{\pm 1})^{s}.(x^-_\beta)^{s+1}&=
 \sum_{m=0}^{s}(x^-_{\beta} t_2^{\pm m})\Lambda^\pm(h_\beta,
\ s-m) +\  X,\\ &(x^+_{\beta}t_2^{\pm
1})^{s+1}.(x^-_\beta)^{s+1}&= \Lambda^\pm(h_\beta, s+1) +Y
 ,\end{eqnarray*} where $X$ and $Y$ are in the left ideal
 of $\butor$ generated by the subalgebra $\gtor(>)$.

 \hfill\qedsymbol\end{lem}

\section{Representation theory of $\gfin$ and $\gaff$}
\noindent In this section, we discuss  the representation theory
of $\gfin$ and $\gaff$. We shall be interested in the
finite--dimensional representations of $\gfin$ and their
analogues, the integrable representations of $\gaff$.

\subsection{}
For $\lambda=\sum_i \lambda(h_i)\omega_i\in P_{fin}^+$, (resp.
$\lambda=\sum_{i=0}^n\lambda(h_i)\omega_i\in P^+_{aff}$), let
$V_{fin}(\lambda)$ (resp. $V_{aff}((\lambda)$) be the unique
irreducible $\gfin$--module (resp. $\gaff$--module) with highest
weight $\lambda$ and highest weight vector $v_\lambda$. Thus,
$V_{fin}(\lambda)=\bufin.v_\lambda$ (resp.
$V_{aff}(\lambda)=\buaff.v_\lambda$) is generated by the element
$v_\lambda$, subject to the following relations:$$\npfin.
v_\lambda =0, \ \ h.v_\lambda=\lambda(h).v_\lambda,\ \
(x^-_{i})^{\lambda(h_i) +1}.v_\lambda =0, \ \ i=1,\cdots ,n,\
h\in\hfin, $$
$$({\text{resp.}}\ \ \npaff. v_\lambda =0, \ \
h.v_\lambda=\lambda(h).v_\lambda,\ \ (x^-_{i})^{\lambda(h_i)
+1}.v_\lambda =0, \ \ i=0,\cdots ,n,\ h\in\heaff.)$$ It is
well--known that the set $\{V_{fin}(\lambda):\lambda\in
P^+_{fin}\}$ is in one--one correspondence with the set of
isomorphism classes of irreducible finite--dimensional
representations of $\frak g_{fin}$. Further, any
finite--dimensional $\gfin$--module is isomorphic to a direct sum
of irreducible $\gfin$--modules.

\subsection{}
 The $\geaff$--modules $V_{aff}(\lambda)$ are not
finite--dimensional, but are integrable \cite{K} in the following
sense.
\begin{defn} A
$\geaff$--module $V$ is said to be integrable if
$$V=\bigoplus_{\mu\in\left(\heaff\right)^*} V_\mu,$$ where
$V_\mu=\{v\in V:hv=\mu(h)v \ \ \forall \ h\in\heaff\}$, and if for
all $\alpha\in R^+_{fin}$, $r\in\bz$ the elements
$x^\pm_{\alpha}t_1^r$ act locally nilpotently on
$V$.\hfill\qedsymbol\end{defn}

To describe the isomorphism classes of irreducible integrable
representations of $\geaff$, we need to introduce two more
families of modules. The first one is obtained by just taking the
restricted dual $V^*_{aff}(\lambda)$ of $V_{aff}(\lambda)$, i.e.,
the $\geaff$--module generated by an element $v^*_{\lambda}$
subject to the relations $$\nmaff. v_\lambda =0, \ \
h.v_\lambda=-\lambda(h).v_\lambda,\ \ (x^+_{i})^{\lambda_i
+1}.v_\lambda =0, \ \ i=0,\cdots ,n,\ \ h\in\heaff.$$ Notice that
the center $c_1$ acts as a positive integer on $V_{aff}(\lambda)$
(if $\lambda\ne 0$)  and
 as a negative integer on $V^*_{aff}(\lambda)$.

The second family that we need are the loop modules which were
introduced in \cite{CP}. For  $k\ge 1$, $\lambda_1,\cdots
,\lambda_k\in P^+_{fin}$, $a_1,\cdots ,a_k\in\bc^*$, $b\in\bc$,
define a $\geaff$--module structure on
\begin{equation*}L(V_{fin}(\lambda_1)\otimes\cdots\otimes
V_{fin}(\lambda_k))= V_{fin}(\lambda_1)\otimes\cdots\otimes
V_{fin}(\lambda_k)\otimes\bc[t,t^{-1}],\end{equation*} as follows:
\begin{eqnarray*} & c_1(v_1\otimes\cdots\otimes v_k)\otimes
t^s& =0,\\ & d_1(v_1\otimes\cdots\otimes v_k)\otimes
t^s&=(s+b)(v_1\otimes\cdots\otimes v_k)\otimes t^s,\\ &
xt_1^r(v_1\otimes\cdots\otimes v_k)\otimes
t^s&=\left(\sum_{i=1}^ka_i^rv_1\otimes\cdots v_{i-1}\otimes
xv_i\otimes v_{i+1}\otimes\cdots\otimes v_k\right)\otimes
t^{r+s},\end{eqnarray*}
 where $\{v_i\in V_{fin}(\lambda_i): 1\le i\le k\}$, $s,r\in\bz$,
 $x\in\gfin$. Denote this module by $V_{aff}(\bolambda,\boa,b)$.
 The
following result was proved in \cite{CP}
\begin{prop}\label{loop} Let $k\ge 1$, $b\in\bc$,  $\lambda_1,\cdots
,\lambda_k\in P^+_{fin}$
and assume that $a_1,\cdots ,a_k$ are distinct non--zero complex
numbers.
\begin{enumerate}
\item[(i)] The $\geaff$--module $V_{aff}(\bolambda,\boa,b)$ is
irreducible iff for every $r\in\bz$ there exists an integer
$m_r\in\bz$ with  $m_r\neq 0\mod r$ such that
$$\sum_{i=1}^ka_i^{m_r}\lambda_i\neq 0.$$
 \item[(ii)] Suppose that  there exists $r>0$ such that
$\sum_{i=1}^ka_i^m\lambda_i=
 0$ for all $m\neq 0\mod r$. Then the  $\gaff$--module
$L(V_{fin}(\lambda_1)\otimes\cdots\otimes V_{fin}(\lambda_k))$ is
completely reducible. The irreducible submodules are generated by
the elements $v_{\lambda_1}\otimes\cdots\otimes
v_{\lambda_k}\otimes t^\ell$, $0\le \ell<r$.
\end{enumerate}\hfill\qedsymbol
\end{prop}

{\bf Remark.} The characters of these modules have been studied in
\cite{Gr}.  In the case of $sl_n$ a crystal model was obtained for
these modules in \cite{Gr1}.

The next result was proved in \cite{C}, \cite{CP}.
\begin{thm} \label{classify} Let $V$ be an irreducible integrable
$\geaff$-module. Assume also that ${\text{dim}}(V_\mu)<\infty$ for all
$\mu\in\heaff$. Let  $\ell\in\bz$ be such that $c_1v=\ell v$ for
all $v\in V$. Then:
\begin{enumerate} \item[(i)]If $\ell>0$ (resp. $\ell<0$)  there exists
$\lambda\in
P_{aff}^+$ such that $V\cong V_{aff}(\lambda)$ (resp. $V\cong
V^*_{aff}(\lambda)$). \item[(ii)] If $\ell=0$, then there exists
$k\ge 1$, $\lambda_1,\cdots ,\lambda_k\in P^+_{fin}$, $b\in\bc$
and distinct non--zero complex numbers $a_1,\cdots ,a_k$ such that
$V$ is isomorphic to either $V_{aff}(\bolambda,\boa,b)$ or to one
of its irreducible submodules as described in Proposition
\ref{loop}.
\end{enumerate}
\end{thm}

{\bf Remark.} In \cite{J} a closely related family of bounded
admissible modules for $\gaff$  were studied and a similar
classification theorem was obtained.

\vskip 12pt

 The category of integrable representations
is in general far from semisimple. However, one has the following
result, \cite{K}.
\begin{thm}\label{compred} Let $V$ be an integrable
$\geaff$--module such  that ${\text{dim}} V_\mu<\infty$ for all
$\mu\in (\heaff)^*$.  Assume  that  there exists $\mu_1,\cdots
,\mu_s\in P^+_{aff}$ such that $$V_\mu\ne 0\implies \mu\in
\mu_i-Q^+_{aff},$$ $$({\text{resp.}} V_\mu\ne 0\implies \mu\in
-\mu_i-Q^-_{aff},)$$ for some $1\le i\le s$. Then, $V$ is
isomorphic to a direct sum of representations of the form
$V_{aff}(\lambda)$ (resp. $V^*_{aff}(\lambda)$) where $\lambda\in
P^+_{aff}$. In particular, if $\lambda,\mu\in P^+_{aff}$, the
tensor product $V_{aff}(\lambda)\otimes V_{aff}(\mu)$ is a
reducible but completely reducible $\geaff$--module.

\hfill\qedsymbol\end{thm}

\subsection{} We give an  easy   example of an integrable representation
$V=\oplus_{\mu\in(\heaff)^*}V_\mu$ (with $\text{dim}\
V_\mu<\infty$ for all $\mu\in(\heaff)^*$) of $\geaff$ which is
indecomposable and reducible. Such representations are studied in
greater detail in \cite{CPweyl}. Assume that $\gfin=sl_2$ and that
$h,x,y$ is the standard basis of $sl_2$. Let $V$ be the free
module over $\bc[t_1,t_1^{-1}]$ of rank 4, with basis
$v_0,v_1,v_2, w_0$. It is not hard to check (using the Chevalley
presentation \cite{K})  that the following formulae define an
action of $\geaff$ on $V$. The center $c_1$ acts trivially, and
for any $v\in V$, $n\in\bz$, $d_1vt^n=nt^nv$. In addition,
\begin{eqnarray*}&h(v_it_1^r)=(2-2i)(v_it^r),\ \ x(v_it^r)=(3-i)(v_{i-1}t^r),\ \
y(v_it^r)=(i+1)(v_{i+1}t^r),\\ & (xt^{-1})(v_0t^r)=0,\ \
(xt^{-1})(v_1t^r)= 2(v_0t^{r-1}),\ \ (xt^{-1})(v_2t^r)=
(v_1+w_0)t^{r-1},\\
 &\gaff.w_0t^r =0,\\ & (y t)(v_0t^r)=(v_1+w_0)t^{r+1},\ \ (y
t)v_1t^{r} =2v_2t^{r+1},\ \  (y t)v_2=0.\\
\end{eqnarray*}It is clear that this module is generated by the element $v_0$
and
that for all $r\in\bz$,  the elements $w_0t^r$ generate a proper
submodule of it.

\subsection{} We shall need the following result, \cite{CPnew}  in our
study of representations of the double affine algebra.
\begin{thm}\label{tenspr} Let $\lambda\in P^+_{aff}$ and let
$\mu_1,\cdots,\mu_k\in P^+_{fin}$. Assume also that $a_1,\cdots,
a_k$ are distinct non--zero complex numbers. The $\geaff$--module
$V_{aff}(\lambda)\otimes V_{aff}(\bomu,\boa)$ is irreducible if
\begin{enumerate}
\item[(i)] $\sum_{i=1}^ka_i\mu_i\ne 0$,\
 \item[(ii)] there exists $\alpha\in R^+_{fin}$
  such that either
 $$(k+1)\lambda(c_1)<(\mu+\lambda)(h_\alpha),$$
 or $$k\lambda(c)<(\mu^*-\lambda)(h_\alpha),$$
where $\mu=\sum_{i=1}^k\mu_i$ and $\mu^*$ is the highest weight of
the $\gfin$--module that is dual to $V_{fin}(\mu)$.
\end{enumerate}\hfill\qedsymbol
\end{thm}
Notice that, in particular, the theorem implies the existence of
irreducible  integrable representations with infinite--dimensional
weight spaces. A partial converse to the theorem was proved in
\cite{CPsup}. More recently, a converse has been proved in
\cite{Ad} in the case of $sl_2$.

\section{Representations of $\gtor$}
The representation theory of $\getor$ is substantially more
complicated than the affine case, and a  classification theorem of
the kind in Theorem \ref{classify} seems much more difficult;
however, see \cite{E1}, \cite{E2}. The study of $\getor$--modules
is closely related to that of  $\gtor$--modules and  for
simplicity we restrict our attention to the representation theory
of $\gtor$. We shall be interested in the category of integrable
representations of $\gtor$ on which  $c_2$ acts trivially, but
$c_1$ acts non--trivially. Since this category is not semisimple,
we shall discuss both irreducible and indecomposable
representations of $\gtor$.

\subsection{} We begin with the definition of integrable
$\gtor$--modules.
\begin{defn} A representation $V$ of $\gtor$  is called integrable
if $$V=\bigoplus_{\lambda\in(\htor)^*}V_{\lambda}$$ where
$$V_{\lambda}=\{v\in V: h.v=\lambda(h)v\ \forall h\in\htor\},$$
and if  the elements $x^\pm_{\alpha}t_1^{m_1}t_2^{m_2}$ act
locally nilpotently on $V$ for all $\alpha\in R^+_{fin}$,
$m_1,m_2\in\bz$. We say that $V$ is admissible if
$\text{dim}V_\lambda <\infty$ for all $\lambda \in(\htor)^*$.
\hfill\qedsymbol
\end{defn}
 It is clear
that if $V$ is integrable, and $0\ne v\in V_\lambda$, then
$\bufin.v$ is a finite--dimensional $\bufin$--module and that
$\bueaff.v$ is an integrable $\geaff$--module.

\subsection{}   Given
$\lambda_1,\cdots ,\lambda_k\in (\heaff)^*$, $a_1,\dots
,a_k\in\bc^*$, define an action of $\getor$ on
$V_{aff}(\lambda_1)\otimes \cdots\otimes V_{aff}(\lambda_k)$ as
follows: $c_2$ acts as zero and for $x\in\gaff$, $m\in\bz$,
$$xt_2^m.(v_1\otimes\cdots\otimes v_k)=\sum_{j=1}^k
a_j^m(v_1\otimes\cdots\otimes v_{j-1}\otimes xv_j\otimes
v_{j+1}\otimes\cdots\otimes v_k),\ \ v_i\in V_{aff}(\lambda_i),
1\le i\le k. $$ Denote this module by $V_{tor}(\bolambda,\boa)$.
Notice that  $$V_{tor}(\bolambda,\boa) \cong
V_{tor}(\bolambda',\boa ')$$ if and only if there exists a
permutation $\sigma$ of $\{1,\cdots ,k\}$ such that
$$\lambda_i=\lambda'_{\sigma(i)},\ \ a_i=a'_{\sigma(i)}.$$
\begin{lem} For all $k\ge 1$, $\bolambda=(\lambda_1,\cdots
,\lambda_k)\in (P^+_{aff})^k$, $\boa=(a_1,\cdots
,a_k)\in(\bc^*)^k$, the module $V_{tor}(\bolambda,\boa)$ is an
integrable admissible $\gtor$--module which  is irreducible iff
$a_r\ne a_s$ for all $1\le r\ne s\le k$. Let $v_{\lambda_r}$ be
the highest weight vector in $V_{aff}(\lambda_r)$. We have
\begin{eqnarray}\label{hw}&\butor(>).v_{\lambda_1}\otimes\cdots\otimes
v_{\lambda_k}
&=0,\\ & \label{eig} \Lambda^\pm(h_i,u
)(v_{\lambda_1}\otimes\cdots\otimes v_{\lambda_k})&=
\prod_{j=1}^k(1-a_j^{\pm 1}
u)^{\lambda_j(h_i)}(v_{\lambda_1}\otimes\cdots\otimes
v_{\lambda_k}),\ \ \forall\  0\le i\le n.\end{eqnarray}
\end{lem}
\begin{pf} Since $x^\pm_{\alpha}t_1^m$ acts locally nilpotently on
$V_{aff}(\lambda_i)$ for all $\alpha\in R^+_{fin}$, it follows
from the definition that $V_{tor}(\bolambda,\boa)$ is integrable.
The module is admissible because the weight spaces of
$V_{aff}(\lambda_i)$ are finite--dimensional. If $a_r\ne a_s$ for
all $r\ne s$, notice that for all $x\in\gaff$, the element
$x\prod_{s\ne r}(t_2-a_s)$ acts only on the $r^{th}$ component of
the tensor product. The irreducibility of
$V_{tor}(\bolambda,\boa)$ now follows from the irreducibility of
the $\geaff$--module  $V_{aff}(\lambda_i)$. Suppose now that
$a_r=a_s$ for some pair $(r,s)$; we can assume without loss of generality (by
applying a suitable permutation) that $r=1,s=2$.  Recall from
\cite{K} that $V_{aff}(\lambda_1)\otimes V_{aff}(\lambda_2)$ is a
reducible $\geaff$--module; let $W$ be a  proper submodule. It is
clear that $W\otimes V_{aff}(\lambda_3)\otimes\cdots\otimes
V_{aff}(\lambda_k)$ is a proper non--zero  submodule of
$V_{tor}(\bolambda,\boa)$.

Finally, note that \eqref{hw} and \eqref{eig} follow from the
definition of the $\gtor$--action on $V_{tor}(\bolambda,\boa)$ and
\eqref{power}
\end{pf}

\subsection{} We now  construct a family of integrable indecomposable
$\gtor$--modules. This family will be maximal in suitable sense.
We begin with the following definition.

Given a collection $\bop=\{p^\pm_r\in(\haff)^*:r\ne 0\}$ and
$\lambda\in (\htor)^*$, let $M_{tor}(\lambda, \bop)$ be the left
$\butor$--module defined as  follows. Let $I(\lambda,\bop)$ be the
left ideal in $\butor$ generated by $$\butor(>),\ \ \ \
\Lambda^\pm(h_i, r)-p^\pm_r(h_i), \ \ \, h'-\lambda(h'),\ \
(h'\in\htor ,0\le i\le n).$$ Set
$$M_{tor}(\lambda,\bop)=\butor/I(\lambda,\bop).$$ Clearly
$M_{tor}(\lambda,\bop)$ is a left $\butor$--module. Let
$v_{\lambda,\bop}$ be the image of $1$ in $M_{tor}(\lambda,\bop)$.
Standard arguments prove that $M_{tor}(\lambda,\bop)$ is a free
$\butor(<)$ module generated by $v_{\lambda,\bop}$ and hence that
$M_{tor}(\lambda,\bop)$ has a unique irreducible quotient, which
we denote by $V_{tor}(\lambda,\bop)$.

The following result determines a necessary and sufficient
condition for $M_{tor}(\lambda,\bop)$ to have an integrable
quotient and identifies the irreducible quotient in this case.

\begin{prop}\label{intquot} The module $M_{tor}(\lambda,\bop)$ has an integrable
quotient iff there exists an $(n+1)$--tuple of polynomials (with
constant term one) $\bpi=(\pi_0,\cdots ,\pi_n)$ in an
indeterminate $u$, such that the following conditions hold for all
$i=0,\cdots ,n$:
\begin{enumerate}
\item[(i)] $\lambda(h_i)={\text{deg}}\ \pi_i$,
  \item[(ii)] $\sum_{r\ge 0}p^\pm_r(h_i)u^r = \pi_i^\pm(u),$
  where $\pi_i^+(u) =\pi_i$ and
  $\pi_i^-(u)=
\frac{u^{\text{deg}\pi_i}\pi^+(u^{-1})}{(u^{\text{deg}\pi_i}\pi^+(u^{-1}))(0)}$.
\end{enumerate}
\end{prop}
\begin{pf} Suppose that $M_{tor}(\lambda,\bop)$ has an integrable
$\gtor$--quotient $W$ and let $w_\bop$ be the image of
$v_{\lambda,\bop}$ in $W$. By the representation theory of $sl_2$
it follows  that $\lambda(h_i)$ is a non--negative integer  for
all $0\le\le n$. Further  $\lambda(h_i)=r_i$ is the  smallest
non--negative integer such that
$$(x^-_i)^{r_i+s}.w_\bop=0, \ \ 0\le i\le n,\ \ \forall s> 0.$$

 Applying
$(x^+_it_2)^{r_i+s}$ to the preceding equation and using Lemma
\ref{deflambda}, we see that $$\Lambda^+({h_i,m}).w_\bop
=p^+_m(h_i)=0\ \ {\text{if}} \ m\ge r_i+1.$$ To see  that
$p^+_{r_i}(h_i)\ne 0$, note that for all $0\le i\le n$ we have,
$$(h_it_2^{-1})(x^-_i)^{r_i+1}.w_\bop=(x_i^-t_2^{-1})(x^-_i)^{r_i}.w_\bop=0.$$
Since the elements $x_i^\pm t_2^{\pm 1}$ generate a subalgebra
isomorphic to $sl_2$ it follows from the representation theory of
$sl_2$ that,$$(x^+_it_2)^{r_i} (x^-_i)^{r_i}.w_\bop\ne 0.$$ Using
Lemma \ref{deflambda} we see that this means that $p^+_r(h_i)\ne
0$ for all $0\le i\le n$. Hence  $\pi^+_i(u)=\sum_s
p^+_{s}(h_i)u^s$ is a polynomial of degree $r_i$. Similarly one
can prove that $\pi^-_i(u)=\sum_sp^-_s(h_i)u^s$ is a polynomial of
degree $r_i$. To see that $\pi_i^\pm$ are related as in the
proposition, one proceeds as in the proof of Proposition
(1.1)(iv),(v) in \cite{CPweyl}, we omit the details.

Conversely, given $\bpi=(\pi_0,\cdots ,\pi_n)$, consider the set
$\{a_1,a_2,\cdots ,a_r\}$ of distinct  roots of
$\pi=\prod_{j=0}^n\pi_j$. Let $m_{ij}$ be the multiplicity with
which $a_i$ occurs as a root in $\pi_j$. Set
$$\mu_j=\sum_{i=0}^nm_{ji}\omega_i, \ \ 1\le j\le r.$$ It is not
hard to see that there exists a map of $\gtor$--modules
$M_{tor}(\lambda,\bop)\to V_{tor}(\bomu,\boa)$. Since
$V_{tor}(\bomu,\boa)$ is integrable, the theorem follows.
\end{pf}

\subsection{} We shall only be interested in the modules
$M_{tor}(\lambda,\bop)$  which are given by an $(n+1)$--tuple of
polynomials as in the preceding proposition. Set
$\lambda_\bpi=\sum_{i=0}^n(\text{deg}\ \pi_i)\omega_i$. Then
$\lambda_\bpi=\lambda$, and we  set
$M_{tor}(\bpi)=M_{tor}(\lambda_\bpi,\bop)$ and denote by
$V_{tor}(\bpi)$ the unique irreducible quotient of
$M_{tor}(\bpi)$. Let $v_\bpi=v_{\lambda,\bop}$. Clearly
$V_{tor}(\bpi)\cong V_{tor}(\bomu,\boa)$ for a suitable choice of
$\bomu,\boa$.

  We now define
the maximal integrable quotient $W_{tor}(\bpi)$ of
$M_{tor}(\bpi)$. Thus, let $W_{tor}(\bpi)$ be the quotient of
$M_{tor}(\bpi)$ by the submodule generated by the elements
$$(x_i^-)^{r_i+1}.v_\bpi,\ \ r_i={\text{deg}}\pi_i,\ \ 0\le i\le
n.$$ Let $w_\bpi$ be the image of $v_\bpi$ in $W_{tor}(\bpi)$. The
following lemma is immediate.
\begin{lem} The $\buaff$--submodule of $W_{tor}(\bpi)$ generated
by $w_\bpi$ is isomorphic to $V_{aff}(\lambda_\bpi)$, where
$\lambda_\bpi=\sum_{i=0}^n ({\text {deg}}\ \pi_i)\omega_i\in
P^+_{aff}.$ \hfill\qedsymbol\end{lem}

\begin{prop}\label{intquot1} The $\gtor$--module $W_{tor}(\bpi)$ is integrable
and
admissible. Further, any integrable quotient of $M_{tor}(\bpi)$ is
a quotient of $W_{tor}(\bpi)$.
\end{prop}
\begin{pf} To see that $W_{tor}(\bpi)$ is integrable, first observe that
$W_{tor}(\bpi)=\butor(<).w_\bpi$. Since the elements $x^\pm_\alpha
t_1^{m_1}t_2^{m_2}$,
 $\alpha\in R^+_{fin}$, $m_1,m_2\in\bz$, act locally nilpotently (via the
adjoint action) on
 $\butor$ it is enough to prove that they
   act locally nilpotently on $w_\bpi$. If $m_2=0$,  the result
   follows since $V_{aff}(\lambda)$ is an integrable $\geaff$--module, and so we
have
   $$(x^-_\alpha t_1^{m_1})^N.w_\bpi=0,$$
   for some $N\ge 0$ depending on $\alpha$ and $m_1$. Applying
   $h_\alpha t_2^{r}$ to the preceding equation, we get
   $$(x^-_\alpha t_1^{m_1})^{N-1}(x^-_\alpha t_1^{m_1}t_2^{r}).w_\bpi=0.$$
Repeating, we find that for any $r_1,\cdots ,r_N\in\bz$
$$(x^-_\alpha t_1^{m_1}t_2^{r_1}) (x^-_\alpha
t_1^{m_1}t_2^{r_2})\cdots (x^-_\alpha
t_1^{m_1}t_2^{r_N}).w_\bpi=0. $$ This proves that  the elements
$(x^-_\alpha t_1^{m_1}t_2^{m_2})$ act nilpotently on $w_\bpi$ and
hence that $W_{tor}(\bpi)$ is integrable.

 To see that $W$ is admissible, fix a total
order on the set $R^+_{aff}$. By the Poincare--Birkhoff--Witt
theorem it follows that  the weight space $W_{\lambda-\eta}$ for
$\eta\in Q^+_{aff}$ is spanned by elements of the form
$$x^-_{\beta_1}t_2^{k_1}\cdots x^-_{\beta_r}t_2^{k_r}.w_\bop,$$
where $r\ge 1$,  $\beta_1\le \beta_2\le\cdots\le\beta_r$ are
elements of $R^+_{aff}$ such that $\sum_s\beta_s=\eta$ and
$k_1,\cdots ,k_r\in\bz$. Since the number of roots in $R^+_{aff}$
that add up to a fixed $\eta\in Q^+_{aff}$ is finite, to see that
$W_{\lambda-\eta}$ is finite--dimensional it suffices to prove
that the number of choices for $k_1,\cdots ,k_r$ is finite.  In
fact, it  suffices to prove that for every $\beta\in R^+_{aff}$
there exists an integer $N_\beta$ such that
$x^-_{\beta}t_2^k.w_\bpi$  is in the span of elements of the
form $x^-_\beta t_2^s.w_\bpi$, with $-N_\beta\le s\le N_\beta$.
For then an obvious induction on $r$ proves that, for every
$\eta\in Q^+_{aff}$, there exists an integer $N_\eta$ such that
$W_{\lambda-\eta}$ is spanned by elements of the form
 $$x^-_{\beta_1}t_2^{k_1}\cdots
x^-_{\beta_r}t_2^{k_r}.w_\bpi,\ \ -N_\eta\le k_i\le N_\eta,
\sum_s\beta_s=\eta. $$

 If $\beta\in R^+_{aff}\backslash
\{m\delta_1:m\in\bz, m>0\}$, it follows from Lemma \ref{deflambda}
and the fact that $(x^-_\beta)^{\lambda(h_\beta)+1}=0$    that
$x^-_\beta t_2^k.w_\bpi$ is in the span of $\{x^-_\beta
t_2^r.w_\bpi: -\lambda(h_\beta)\le r\le \lambda(h_\beta)\}$. If
$\beta=m\delta_1$, $m>0$, the corresponding negative root vectors
in $\geaff$ are $\{h_it_1^{-m}:1\le i\le n\}$. Take
$\beta=\alpha_i+m\delta_1$ in   the first equation in Lemma
\ref{deflambda} and apply $x_i^+$ to it. This gives,
$$\sum_{r=0}^{s}h_it_1^{-m}t_2^s\Lambda^\pm(h_i, s-r).w_\bpi =0.$$
Again, it follows that the element $h_it_1^{-m}t_2^{m_2}.w_\bpi$
is in the span of the elements $\{h_it_1^{-m}t_2^{s_2}.w_\bpi:
-\lambda(h_i)\le s_2\le \lambda(h_i)\}$. This completes the proof
of the proposition.

\end{pf}
\begin{cor}\label{gplus}  We have $W_{tor}(\bpi)=\bu(\geaff\otimes
\bc[t_2]).w_\bpi$.
\end{cor}
\begin{pf} Since $W_{tor}(\bpi)$ is an integrable module, it
follows from Lemma \ref{deflambda} that if $\beta\in
R^+_{aff}\backslash \{r_1\delta:r_1\in\bz\}$, then
$$\sum_{m=0}^{N_\beta}(x^-_{\beta} t_2^{ m})\Lambda^\pm(h_\beta, \
\ N_\beta-m).w_\bpi=0,$$ where $N_\beta=\lambda_\bpi(h_\beta)$.
Applying $h_\beta t_2^{m_2}$, $m_2<0$, to both sides of the
equation gives $$\sum_{m=0}^{N_\beta}(x^-_{\beta} t_2^{
m+m_2})\Lambda^\pm(h_\beta, \ N_\beta-m).w_\bpi=0.$$ This proves
that  $(x^-_\beta t_2^{m_2}).w_\bpi$ is in the span of elements
$\{(x^-_\beta t_2^{m}).w_\bpi:m>m_2\}$ and hence by a simple
induction in the span of elements of the form  $\{(x^-_\beta
t_2^{m}).w_\bpi:m\ge 0\}$. A similar statement then follows for
roots of the form $\{r_1\delta_1:r_1\in\bz\}$, exactly as in the
proof of the proposition. The corollary follows.
\end{pf}



We conclude this section with some  result on the structure of the
modules $W_{tor}(\bpi)$ that we shall need in the last
section of the paper.

\begin{lem}  Regarded as  a  module  for $\gaff$, we have
 $$W_{tor}(\bpi)\cong\bigoplus_{\mu\in
P^+_{aff}} m(\mu) V_{aff}(\mu),$$ for some $m(\mu)\in \bz$,
$m(\mu)\ge 0$. Further, $m(\lambda_\bpi)=1$ and $$m(\mu)\ne
0\implies \mu=\lambda_\bpi-P^+_{aff}.$$\end{lem}
\begin{pf} To prove (i), notice that
$$\left(W_{tor}(\bpi) \right)_\mu\ne 0\implies
\mu\in\lambda_\bpi-Q^+_{aff},$$ and
${\text{dim}}\left(W_{tor}(\bpi)\right)_\mu$ is finite. The Lemma
now follows from Proposition \ref{compred}.\end{pf}

Although, we have so far been interested only in the
representations of $\gtor$, we shall also need the corresponding
results for the algebra $\gcurr=\gaff[t_2]$ of polynomial maps
$\bc\to\gaff$. One can define, in the obvious way, a family of
$\gcurr$--modules $M_{tor}^+(\lambda,\bop^+)$, where
$\lambda\in(\heaff)^*$, $\bop^+ =\{p^+_r\in(\haff)^*:r\in\bz,r
>0\}$. It is not hard to see that there exists an injective map of
$\gcurr$--modules $$M^+_{tor}(\lambda,\bop^+)\to
M(\lambda,\bop),$$ where $\bop=\{p_r:r\in\bz, r\ne 0\}$ is such
that $p_r=p_r^+$ for all $r>0$. Further, one can show as in
Proposition \ref{intquot}, that $M^+_{tor}(\lambda,\bop^+)$ has an
integrable quotient if and only if $\bop^+$ is  given by an
$(n+1)$--tuples of polynomials. Let $W^+_{tor}(\bpi)$ be the
corresponding maximal integrable quotient, i.e. the quotient of
$M^+_{tor}(\bpi)$ by the submodule generated by elements of the
form $(x^-_i)^{\text{deg}\ \pi_i +1}.v_\bpi$.
\begin{prop}\label{gcur} Assume that $\bpi$ is an $(n+1)$--tuple of
polynomials with constant term one. Then, $$W^+_{tor}(\bpi)\cong
W_{tor}(\bpi)$$ as $\gcurr$--modules.\end{prop}
\begin{pf} It is clear that there exists a map of $\gcurr$--modules
$\iota:M^+_{tor}(\bpi)\to W_{tor}(\bpi)$, which by Corollary
\ref{gplus} is surjective. Suppose that $v=gv_\bpi\in
M^+_{tor}(\bpi)$ maps to zero. This means that
$\iota(gv_\bpi)=g\iota(v_\bpi)$ is in the $\gtor$--submodule of
$M_{tor}(\bpi)$ generated by the elements $(x^-_i)^{\text{deg}\
\pi_i +1}.\iota(v_\bpi)$. Since
$$\sum_{i=0}^n\bu(\gtor(<)).(x^-_i)^{\text{deg}\ \pi_i
+1}.\iota(v_\bpi)=\sum_{i=0}^n\bu(\nmaff\otimes\bc[t_2^{-1}]t_2^{-1})\bu(\nmaff\otimes\bc[t_2])(x^-_i)^{\text{deg}\
\pi_i +1}\iota(v_\bpi),$$ and $M_{tor}(\bpi)$ is  a free
$\gtor(<)$--module, it follows from the  PBW--theorem that
$$g\in\sum_{i=0}^n\bu(\nmaff\otimes\bc[t_2])(x^-_i)^{\text{deg}\ \pi_i
+1}\iota(v_\bpi).$$ Hence, we find that the induced map
$W_{tor}^+(\bpi)\to W_{tor}(\bpi)$ is injective, and the
proposition is proved.

\end{pf}

The final result of this section  is an analog for the modules
$W_{tor}(\bpi)$ of  the factorization that holds for the
irreducible modules $V_{tor}(\bomu,\boa)$. The proof of this
Proposition is a modification of the proof of Proposition 3.1 in
\cite{CPweyl}, the details of the proof can be found in
\cite{thesis}.

\begin{prop}\label{factor}
 Let $a_1,\cdots ,a_k $ be
the distinct roots of $\prod_{i=0}^n\pi_i$. For $1\le s\le k$, and
$0\le i\le n$, assume that $a_s$ occurs as a root of $\pi_i$ with
multiplicity $m_{i,s}$.
 For $1\le j\le
k$, let $\bpi_j=(\pi_{j,0},\cdots ,\pi_{j,n})$ be defined by,
$$\pi_{j,i}=(1-a_j^{-1}u)^{m_{i,j}}.$$ Then we have an isomorphism
$$W(\bpi)\cong W(\bpi_1)\otimes\cdots \otimes W(\bpi_k),$$ of
$\gtor$--modules.
\end{prop}

\section{An irreducibility criterion for $W_{tor}(\bpi)$}
It follows from Proposition \ref{intquot1} that the unique
irreducible quotient of $W_{tor}(\bpi)$ is $V_{tor}(\bpi)$. In
this section we give a condition for the quotient map to be an
isomorphism.

Recall that $\theta$ is the highest root of $R^+_{fin}$ and write
$h_\theta=\sum_{i=1}^n m_ih_i$. For $0\le i\le n$, and
$a\in\bc^\times$, define an $(n+1)$--tuple of  polynomials
$\bpi_{i,a}$ by
$$\pi_j=1, \ \ j\ne
i,\ \ \pi_i=(1-au).$$
\begin{thm}\label{irred}  Assume that  either  $i=0$ or that $1\le i\le n$ is
such that $m_i=1$. Then, $$W_{tor}(\bpi_{i,a})\cong
V_{tor}(\bpi_{i,a})$$ as $\buaff$--modules.\end{thm}
\begin{pf} Set $\bpi=\bpi_{i,a}$. Consider the map
$$\phi_r: (\gaff t_2^{r}])\otimes  W_{tor}(\bpi)\to
W_{tor}(\bpi)$$ given by $\phi_r(xt_2^r, w) =(xt_2^r)w$. This is
clearly a map of $\gaff$--modules, where $\gaff$ acts on the first
factor through the adjoint representation. Let $p_\mu$ denote the
projection of $$p_\mu:W_{tor}(\bpi)=\bigoplus_{\nu\in P^+_{aff}}
V_{aff}(\nu)\to W(\mu)$$ where $W(\mu)$ is the  $\mu^{th}$
isotypical component of $W_{tor}(\bpi)$. Then $p_\mu$ is a map of
$\gaff$--modules, and hence for every $r\in\bz$ we get  a map of
$\gaff$--modules $$\phi_{r,\mu}=p_\mu.\phi_r:(\gaff t_2^r)\otimes
W_{tor}(\bpi)\to W(\mu).$$  We show that for all $r\in\bz$ the
restriction of $\phi_{r,\mu}$ to $(\gaff t_2^r)\otimes
V_{aff}(\bpi)$ is zero if $\mu\ne \lambda_\bpi$.  Observe that
$\lambda_\bpi=\omega_i$ where $i$ is such that $m_i=1$ and hence
$\omega_i(c_1) =1$. This implies that
$$(\lambda_\bpi+\theta)(h_\theta)>2,\ \ {\text{if}}\\ i\ne 0$$ and
$$(\theta-\lambda_\bpi)(h_\theta)>1, \ \ i=0.$$ Further, since
$c_1 t_2^r.w_\bpi$ is a scalar multiple of $w_\bpi$, it follows
that   $\phi_{r,\mu}(c_1 t_2^r,v)=0$ for all $v\in
V_{aff}(\lambda_\bpi)$ if $\mu\ne \lambda_\bpi$. This implies that
$\phi_{r,\mu}$ factors through to  a map of $\gaff$--modules
$$\gfin( \bc[t_1,t_1^{-1}])t_2^r\otimes V_{aff}(\lambda_\bpi)\to
W(\mu).$$ Since $$\gfin ( \bc[t_1,t_1^{-1}])t_2^r \cong
L(V_{fin}(\theta))$$ as $\gaff$--modules, it follows from Theorem
\ref{tenspr}  that $(\gaff t_2^r)\otimes V_{aff}(\lambda_\bpi)$ is
irreducible. On the other hand since this module has
infinite--dimensional weight spaces, it follows from Schur's lemma
that  $\phi_{r,\mu}=0$. This means that as a $\geaff$--module we
have $$W_{tor}(\bpi_{i,a})\cong V_{aff}(\lambda_i).$$ On the other
hand the irreducible quotient of $W_{tor}(\bpi_{i,a})$ is the
module $V_{aff}(\lambda_i,a)$ and this by definition is isomorphic
to  $V_{aff}(\lambda_i)$  as a $\geaff$--module.  The theorem now
follows.

\end{pf}

\section{Fusion product and reducibility of $W_{tor}(\bpi)$}
To understand the $\gaff$--module structure of $W_{tor}(\bpi)$ and
in particular to prove that it is not isomorphic to
$V_{tor}(\bpi)$, we need to introduce the notion of the fusion
product of representations of $\gaff\otimes\bc[t_2]$. This was
introduced by Feigin and Loktev in \cite{FL}.

Thus, let $\frak a$ be any Lie algebra and let $\frak a[t]$ be the
algebra of polynomial maps $\bc\to\frak a$. For $k\ge 1$, let
$V_1, \cdots , V_k$ be representations of $\frak a$ and let
$a_1,\cdots ,a_k\in\bc$ be arbitrary complex numbers. As in the
previous sections, one sees that  the tensor product
$V_1\otimes\cdots\otimes V_k$ admits a structure of $\frak
a[t]$--module given by $$xt^r(v_1\otimes\cdots \otimes
v_k)=\sum_{j=1}^ka_j^rv_1\otimes\cdots\otimes v_{j-1}
xv_j\otimes\otimes v_{j+1}\cdots\otimes v_k,$$ for all $x\in\frak
a$, $r\ge 0$, $v_i\in V_i$, $1\le i\le k$. Denote this module by
$V_1(a_1)\otimes\cdots\otimes V_k(a_k)$. The following lemma is
easily proved.
\begin{lem} Assume that there exist elements $v_j\in V_j$ such
that $V_j=\bu(\frak a[t]).v_j$ for all $1\le j\le k$ and that
$a_1,\cdots ,a_k$ are distinct complex numbers. Then,
$V_1(a_1)\otimes \cdots\otimes V_k(a_k)$ is generated as an $\frak
a[t]$--module by $v_1\otimes\cdots\otimes v_k$.
\end{lem}

The Lie algebra $\frak a[t]$  and its enveloping algebra
$\bu(\frak a[t])$ are obviously graded by powers of $t$. Let
$\bu(\frak a[t])_r$ be the $r^{th}$--graded piece. If $V_1,\cdots
V_k$ are generated by elements $v_1,\cdots ,v_k$ and $a_1,\cdots,
a_k$ are distinct complex numbers, the $\frak a[t]$--module
$V=V_1(a_1)\otimes \cdots\otimes V_k(a_k)$ admits an $\frak a$ (but
not an $\frak a[t]$) equivariant filtration. The
$r^{th}$--filtered piece is given by $$ V_r=\bigoplus_{0\le s\le
r}\left(\bu(\frak a[t])\right)_s.(v_1\otimes\cdots\otimes v_k).$$
The associated graded vector space $\oplus_{r\ge 0} V_r/V_{r-1}$
is obviously a $\frak a$--module. Since $\frak at^s. V_r\subset
V_{r+s}$, one can define an $\frak a[t]$--module structure on
$\oplus_{r\ge 0} V_r/V_{r-1}$. This module,  denoted by
$V_1(a_1)*V_2(a_2)*\cdots *V_k(a_k)$,  is called the fusion product
of the modules $V_1,\cdots ,V_k$ with respect to $a_1, \cdots
,a_k$. Let $v_1 * \cdots * v_k$ be the image of
$v_1\otimes\cdots\otimes v_k$ in $V_1(a_1)*V_2(a_2)*\cdots
*V_k(a_k)$.

\vskip 12pt

 {\bf Remark} As an $\frak a$--module, it is clear that
\begin{equation}\label{amod} V_1(a_1)*V_2(a_2)*\cdots
*V_k(a_k)\cong V_1\otimes\cdots\otimes V_k.\end{equation}

\vskip 12pt

 Assume now that $\frak a=\gaff$ and that
$V_j=V_{aff}(\lambda_j)$, where $\lambda_j\in P^+_{aff}$, $1\le
j\le k$. The next lemma is not hard to check.

\begin{lem}\label{fus}  Let $k\ge 1$, and let $\lambda_1,\cdots ,\lambda_k\in
P^+_{aff}$.
  Assume
that $a_1,\cdots ,a_k$ are distinct scalars. The fusion product
$V_{aff}(\lambda_1)( a_1)*\cdots *V_{aff}(\lambda_k)(a_k)$ is
generated as an $\frak a[t]$--module by the element $v=v_1
\otimes\cdots \otimes v_k$. Further, $$\npaff[t].v=0,\ \
(\haff)t\bc[t].v=0,\ \ h.v=\sum_j\lambda_j(h).v,\ \
h\in\haff.$$\hfill\qedsymbol
\end{lem}

For any $\lambda=\sum_{i=0}^n \lambda(h_i)\omega_i$,
$\lambda(h_i)\ge 0$, fix distinct complex numbers $\{c_{ij}:0\le
i\le n, 1\le j\le \lambda(h_i)\}$. Consider the fusion product
$$W(\lambda) = V_{aff}(\omega_0 )(c_{01}) *\cdots
*V_{aff}(\omega_0)(c_{0\lambda_0})*\cdots
*V_{aff}(\omega_n)(c_{n1})*\cdots * \cdots
V_{aff}(\omega_n)(c_{n\lambda_n}).$$ For any complex number
$a\in\bc$, let $W(\lambda, a)$ be the $\gtor[t]$--module obtained
by pulling back $W(\lambda)$ through the Lie algebra homomorphism
$\gaff[t]\to\gaff[t]$ given by sending $xt^r\to x(t-a)^r$.

We turn now to the $\gtor$-module
$W_{tor}(\bpi)$. We  restrict our attention to the case
where there exists $a\in\bc^*$ and $n_j\in \bz$, $n_j\ge 0$, such
that if, $\bpi=(\pi_1,\cdots ,\pi_n)$ then,
 $$\pi_j(u)=(1-au)^{n_j},\ \
1\le j\le n.$$

\begin{thm} There exist a surjective map of $\gaff[t_2]$--modules
$W_{tor}(\bpi)\to W(\lambda_\bpi, a)$.
\end{thm}
\begin{pf} By Proposition \ref{gcur} it is enough to prove that
the element $(v_{\omega_0})^{\otimes(\lambda_0)}\otimes \cdots
\otimes (v_{\omega_n})^{\otimes \lambda_n}$ satisfies the defining
relations of the element $w_\bpi$. But this is clear from Lemma
\ref{fus}.
\end{pf}

It is now not hard to see that  $W_{tor}(\bpi)$ is generally
reducible. This is because we know from the results of Section 3
that the irreducible quotient of $W_{tor}(\bpi)$ in this case is
$V_{tor}(\lambda_\bpi,a)$, and $$V_{tor}(\lambda_\bpi,a)\cong
V_{aff}(\lambda_\bpi)$$ as $\geaff$--modules. On the other hand
the module $W(\lambda_\bpi,a)$ is isomorphic as a $\geaff$--module
to a tensor product of integrable irreducible $\geaff$--modules,
and hence must be reducible as a $\geaff$--module.

\end{document}